\documentclass[10pt]{article}
\usepackage{amsmath, amssymb}
\usepackage{amscd}

\newtheorem{theor}{Theorem}[section]

\newtheorem{lem}[theor]{Lemma}
\newtheorem{propo}[theor]{Proposition}
\newtheorem{coro}[theor]{Corollary}

\newenvironment{pf}{{\it Proof.}}{\hfill $\square$\\}

\begin{document}
\title{\bf Leafwise Holomorphic Functions}
\author{R. Feres and A. Zeghib}
\date{June 20, 2001}  
\maketitle
      
\footnotetext{Dept. of Mathematics - 1146, Washington University,
 St. Louis, MO 63130, USA.}
 \footnotetext{UMPA - \'Ecole Normale Sup\'erieure de Lyon,
69364 Lyon CEDEX 07, France.}      
\footnotetext{{\em Mathematical Subject Classification}: Primary
37C85; Secondary 32A99}
\footnotetext{{\em Key words and phrases}: foliated spaces,
leafwise holomorphic functions.}

\begin{abstract}
It is a    well-known and elementary fact 
   that a holomorphic
function      on a compact  complex manifold   without boundary
 is necessarily constant. The  purpose of the present
article  is to investigate whether, or to what extent, a
similar property 
  holds in the setting of  holomorphically foliated
spaces. 
\end{abstract}

\section{Introduction and  Statement of Results}
Suppose that $M$ is   a
 compact manifold,  $\mathcal{F}$ is a continuous
foliation of $M$ by
 (not necessarily compact) complex leaves, 
 and that 
 $f$ is a continuous 
leafwise holomorphic function. 
The question we wish to study is  whether or not
 $f$ must be leafwise constant.

We will actually work in the setting
of {\em foliated spaces}, as defined in 
\cite{candelconlon}. Thus $M$ is only a  topological
space, while the leaves of $\mathcal{F}$ admit a    smooth
manifold structure that varies continuously on $M$.
 The term {\em foliated manifold}
will be reserved for   when $M$ has
a differentiable structure relative to
which each leaf of $\mathcal{F}$ is $C^1$ immersed
and the foliation tangent bundle,
$T\mathcal{F}$, is a $C^0$ subbundle of $TM$. In all
situations, it will be assumed that $M$ is compact and
connected. We also assume that $(M,\mathcal{F})$ is
a {\em holomorphically foliated} space, by which we mean
that each leaf of $\mathcal{F}$  carries the structure
of a complex manifold and that this
structure varies continuously on $M$. 
The foliated space   $(M,\mathcal{F})$ may, on
occasion, also carry   a {\em
leafwise Hermitian} metric (resp., is {\em leafwise
K\"ahler}); that is, the leaves of $\mathcal{F}$ may carry
a smooth Hermitian metric (resp., a K\"ahler metric) that,
together with all its derivatives along leaves,
varies continuously on
$M$.

If the foliation is such that any
continuous leafwise holomorphic function is leafwise
constant  we will say that $(M,\mathcal{F})$ -- or
simply 
$\mathcal{F}$ --  is   {\em holomorphically plain}. It
will be
proved that a number of general  classes of holomorphically
foliated manifolds are holomorphically plain. 
We also give an example of a real analytic, non-plain, 
holomorphically foliated
manifold (with  real analytic leafwise holomorphic
functions that are not leafwise constant) and indicate a
general construction that shows how such examples can be obtained.

When the leaves of $\mathcal{F}$, individually,
 do not  
support non-constant bounded holomorphic   functions, then clearly
$\mathcal{F}$ is holomorphically plain. This is the case,  for
example, if the universal covering spaces  of the leaves 
are isomorphic as complex manifolds to $\mathbb{C}^n$.
  Therefore, the
question only becomes meaningful for cases such as,
say, a  foliation 
 by Riemann surfaces
of hyperbolic type, for which the leaves do support
non-constant bounded holomorphic
functions. It is then necessary
to understand the constraining
role played by the foliation dynamics.

The subject of this article has a  natural  
counterpart for foliations with leafwise
Riemannian metrics and functions that are leafwise
harmonic. There is, as well, a discretized form
of the   problem (of deciding whether leafwise harmonic
functions are  leafwise constant) in the
setting of actions of finitely generated groups and
functions that are harmonic along orbits for a
combinatorial Laplacian. These ``harmonic'' variants
of the subject contain some essential additional
difficulties and, for the sake of keeping this
article as elementary as possible, they will be treated
elsewhere.

The first author wishes to thank the members of the UMPA - ENS
Lyon,  and FIM - ETH Z\"urich for their hospitality
while this work was being written.

\subsection{Group Actions and Foliated Bundles}
Holomorphically foliated spaces   arise in a number
of ways. For example, as the orbit foliation of
a   locally free $C^1$ action of a complex connected
Lie group on a compact manifold, or as a {  foliated
bundle} over a compact connected complex manifold.
See also the work of E. Ghys \cite{ghys2}
on laminations by Riemann surfaces and leafwise meromorphic
functions.

It is a rather easy fact that if $(M, \mathcal{F})$ is
the orbit
foliation of a continuous locally free action of  a connected 
complex Lie group $G$, then $\mathcal{F}$ is holomorphically
plain, where the complex structure on leaves is the one that makes
the orbit map $o_x:g \mapsto gx $ (from $G$ onto the leaf
of $x$) a local isomorphism of complex manifolds. 
We have  the following slightly more general fact.
\begin{propo}\label{easy} Let $(M,\mathcal{F})$ be a
holomorphically
  foliated space  such that $T\mathcal{F}$ is
 holomorphically trivial. Then $\mathcal{F}$ is
holomorphically plain.
\end{propo}
  
The next corollary follows from the proposition
and a simple remark due to Wang which is explained later.
 
\begin{coro}\label{coroeasy}
If $(M,\mathcal{F})$ is a holomorphically foliated
space having a dense leaf and $T\mathcal{F}$ is
holomorphically trivial, then $\mathcal{F}$ is the orbit
foliation of a   locally free action
of a connected complex Lie group. 
\end{coro}

A much more interesting class of examples 
 consists of
foliated bundles. We now recall the definition of
   foliated bundles in the
special setting that concerns us here. Let
$S$ denote a compact connected complex manifold, let
$\tilde{S}$ be the universal covering space of $S$, and denote 
by $\Gamma$ the group of deck transformations of $\tilde{S}$.
Let $X$ be a compact connected space on which $\Gamma$ acts
by homeomorphisms. The action can be represented by
a   homomorphism $\rho:\Gamma\rightarrow \text{Homeo}(X)$ from
$\Gamma$ into the group of homeomorphisms of $X$.
Then $\Gamma$ acts on the product $\tilde{S}\times X$ 
in the following way:
$(p,x)\gamma:=(p\gamma, \gamma^{-1}(x))$, for $p\in \tilde{S},
x\in X$ and $\gamma\in \Gamma$. Let $M:=(\tilde{S}\times X)/\Gamma$
be the space of $\Gamma$-orbits. The natural
projection $\pi:M\rightarrow S$ gives $M$ the structure of
a fiber bundle over $S$ whose fibers are homeomorphic to $X$,
and $M$ is foliated by complex manifolds, transversal to
the fibers of $\pi$, which are coverings of
$S$. The resulting foliated space will be written
$(M_\rho,\mathcal{F}_\rho)$.

Note that   the dynamics of a foliated
bundle is largely determined by the dynamics of the
$\Gamma$-action on $X$. Thus, for example, if $X$ admits a
$\Gamma$-invariant finite measure of full support (or
if the $\Gamma$-action has a unique minimal set), then
$M_\rho$ admits a completely invariant finite measure of
full support (resp.,  $(M_\rho,\mathcal{F}_\rho)$ has
a unique minimal set). It should also be clear that
  $x\in 
X$ has a finite $\Gamma$-orbit if and only if
$\tilde{S}\times\{x\}$ maps to a closed leaf of
$\mathcal{F}_\rho$.

\subsection{Leafwise Hermitian and K\"ahler Foliations}
The results in this section assume that 
$(M,\mathcal{F})$ is provided with a (continuous on $M$) 
leafwise Hermitian metric  
$\langle\cdot,\cdot\rangle$. We
denote by $\Omega$ the associated leafwise volume form
and define the divergence of a
continuous
leafwise smooth vector field $X$ to be the function 
$\text{div}X$ such that 
  $\mathcal{L}_X\Omega=(\text{div}X) \Omega$,
where $\mathcal{L}_X$ denotes the Lie derivative along $X$.
A   Borel measure $m$ on $M$ is said to be 
{\em completely invariant} if  
$$m(\text{div}X):=\int_M \text{div}X(x)\, dm(x)=0$$
for all continuous, leafwise 
smooth   vector field $X$.
Completely invariant measures 
are equivalent to {\em holonomy invariant} transverse
measures on $(M,\mathcal{F})$ (cf. \cite{connes}).

A measure $m$ on $M$ will be said to have {\em full support}
if its support coincides with $M$.

\begin{propo}\label{corototransmeas} 
If $(M,\mathcal{F})$ is a   leafwise K\"ahler
 foliated space that admits a
completely invariant measure of full support, then
$\mathcal{F}$ is holomorphically plain.
\end{propo}

 Let
$\bar{\partial}^*$ denote the adjoint operator to
$\bar{\partial}$ with respect to the chosen leafwise
Hermitian metric, and define 
the $\bar{\partial}$-Laplacian on leafwise smooth
functions by
$\Delta_{\bar{\partial}}=\bar{\partial}^*\bar{\partial}$
(we use the notations and sign conventions of
\cite{kodaira}). A Borel measure
$m$ on
$M$ will be called {\em
$\Delta_{\bar{\partial}}$-harmonic}, or simply {\em
harmonic}, if 
$$ m(\Delta_{\bar{\partial}}h):=\int_M
(\Delta_{\bar{\partial}}h)(x)\, dm(x)=0$$ for all continuous,
leafwise smooth function $h:M\rightarrow\mathbb{C}$.

Proposition \ref{corototransmeas} is an immediate
consequence of the next proposition, itself an
immediate consequence of 
\cite[Theorem 1(b)]{garnett}. (We remark that
on a K\"ahler manifold
$\Delta_{\bar{\partial}}h=-\frac12\text{div}\,\text{grad} h$.)

\begin{propo}\label{transmeas}  
Suppose that $(M,\mathcal{F})$  is a   leafwise
Hermitian foliated
space. Also suppose
that   the union of the supports of  all
$\Delta_{\bar{\partial}}$-harmonic measures is $M$. 
 Then   $(M,\mathcal{F})$ 
  is
holomorphically plain.
\end{propo}

\subsection{Foliations with Few Minimal Sets}

Unless  
specified otherwise, $(M,\mathcal{F})$ will continue to denote a 
compact, connected,   holomorphically foliated space.
We recall that  a {\em minimal set} $X$  of
$(M,\mathcal{F})$  is a closed, non-empty,
$\mathcal{F}$-saturated subset of $M$ that has no proper
subset with these same properties.
 If $M$ is, itself, a minimal
set, then $\mathcal{F}$ is said to be a {\em minimal
foliation}.

\begin{propo}\label{prop1}
Suppose that the closure of each
leaf of $(M,\mathcal{F})$ contains
(at most) countably many   minimal sets.
Then $\mathcal{F}$ is holomorphically plain.
\end{propo}

Proposition \ref{prop1} clearly applies  to minimal
foliations. Also note that if  $\mathcal{F}$ is a
(transversely) Riemannian foliation, then the closure of
each leaf is a minimal set \cite{molino}, so the
proposition also applies. 
Therefore one has the next corollary, which will be used later a
number of times. (Note that a compact group of diffeomorphisms
of a compact manifold
must preserve a Riemannian metric.)
\begin{coro}\label{transmeaslemma}
Let $(M_\rho, \mathcal{F}_\rho)$ be a holomorphically
foliated   bundle over   $S$ with fiber a compact manifold
$X$,
where the homomorphism $\rho$ has values in a compact
group of diffeomorphisms  of  $X$.
Then
$\mathcal{F}_\rho$ is holomorphically plain.
\end{coro}

There is also a large class of foliated bundles
associated to $\Gamma$-actions on projective space
$\mathbb{F}P^n$, $\mathbb{F}=\mathbb{R}$ or $\mathbb{C}$,
and derived from linear representations of $\Gamma$,
for which the hypothesis of the proposition are satisfied.
 This will be described after introducing some notations.
Let
$GL(W)$ be  the group of linear automorphisms
of $W$,
where $W$ is a vector space over $\mathbb{F}$.   
The quotient of $GL(W)$ by its center will be 
written $PGL(W)$ and
the projective space associated to $W$ will
be written $P(W)$. If
$W=\mathbb{F}^n$, we write $PGL(W)=PGL(n,\mathbb{F})$
and $P(W)=\mathbb{F}P^{n-1}$.

 Let $\rho:\Gamma\rightarrow
GL(W)$ be a linear representation of a
group
$\Gamma$ on a vector space $W$. An element $\gamma\in
\Gamma$ will be called {\em proximal} if the maximal
characteristic exponent of $\rho(\gamma)$ is simple.
The next result is a consequence of Proposition
\ref{convergence} and \cite[3.4 and 3.6,
Chapter VI]{margulis}.

\begin{propo}\label{zariski}
Let $S$ be a connected, compact, complex manifold
with fundamental group $\Gamma$, $W$ an $n$-dimensional
vector space over $\mathbb{F}$,
and $\rho:\Gamma\rightarrow GL(W)$ a
continuous homomorphism for which $\Gamma$ contains a
proximal element.    Let $(M_\rho,\mathcal{F}_\rho)$ be
the  foliated bundle over 
$S$ with fiber $P(W)$, where $\Gamma$ acts
on $P(W)$ via $\rho$. Then 
$(M_\rho,\mathcal{F}_\rho)$ is  holomorphically plain.
\end{propo}

The hypothesis of
Proposition \ref{zariski}   holds if the image of
$\Gamma$ in $PGL(W)$ is Zariski dense and not precompact
\cite[Theorem 4.3(i)]{margulis}. If the image is precompact,
we can apply Corollary \ref{transmeaslemma}, so the 
following corollary holds.   

\begin{coro}\label{zariskidense}
Let $S$ be a connected, compact, complex manifold
with fundamental group $\Gamma$, $W$ an $n$-dimensional
vector space over $\mathbb{F}$,
and $\rho:\Gamma\rightarrow GL(W)$ a continuous
homomorphism such that the image of $\Gamma$ 
in $PGL(W)$ is Zariski dense.
 Let $(M_\rho,\mathcal{F}_\rho)$ be
the  foliated bundle over 
$S$ with fiber $P(W)$, where $\Gamma$ acts
on $P(W)$ via $\rho$. Then 
$(M_\rho,\mathcal{F}_\rho)$ is  holomorphically plain.
\end{coro}

It will   follow from Corollary \ref{zariskidense} 
that  foliated bundles    
 associated to  projective (linear) actions
of $\pi_1(S)$, for a compact Riemann surface
$S$, are generically holomorphically plain, in the sense
described below.
 
We first recall some definitions.
Let $S=\mathbb{D}/\Gamma$ be a surface of genus $g\geq 2$,
where $\mathbb{D}$ denotes the Poincar\'e disc and
$\Gamma$ is a cocompact discrete group (without torsion) of
hyperbolic
 isometries. 
Let $G$ be  an algebraic group and denote by
$\text{Hom}(\Gamma,G)$ the variety of    homomorphisms
from
$\Gamma$ to $G$. The structure of algebraic variety is 
obtained by identifying 
$\text{Hom}(\Gamma,G)$ with a   subvariety of $G^{2g}$
defined by equations representing relations among elements
in a generating set for $\Gamma$.

\begin{theor}\label{generic} Let $G=GL(n,\mathbb{C})$.
Then there
  is a Zariski open dense   subset
$U$ in $\text{\em Hom}(\Gamma, G)$ such
that,    for each $\rho\in U$,
  the foliated bundle
 $(M_\rho,\mathcal{F}_\rho)$ 
for the  corresponding $\Gamma$-action on
$\mathbb{C}P^{n-1}$  is
holomorphically plain.
\end{theor}

The key step in the results mentioned in this section
is the proposition given next.
 A topological space $X$ equipped with an action
of a group $\Gamma$ by homeomorphisms will be called here
a {\em convergence $\Gamma$-space}, or a $\Gamma$-space
of {\em convergence type}, if the following holds:
there exists a (at most) countable family of
subsets $X_i\subset X$, $i=1,2,\dots$, such that
(i) the intersection of all the $X_i$ is  (at most)
countable   and (ii) for each $i$  there is a sequence
$\gamma_m\in \Gamma$ and a point $x_i\in X$ such
that  $\gamma_m(y)$ converges to $x_i$ as $m\rightarrow
\infty$, for each $y$ in the complement of $X_i$.

As a simple example, let 
$X=\mathbb{C}P^1$ and 
$\rho:\Gamma\rightarrow
PSL(2,\mathbb{C})$ a homomorphism such that $\rho(\Gamma)$
is not relatively compact. Then $X$, with the 
$\Gamma$-action obtained from $\rho$, is a convergence
$\Gamma$-space. The $\Gamma$-actions on $\mathbb{F}P^{n-1}$
of Proposition \ref{zariski} 
as well as the natural action of any unbounded subgroup
$\Gamma$ of a Gromov-hyperbolic group
$G$ on the   boundary $\partial G$, also  define
convergence
$\Gamma$-spaces.

\begin{propo}\label{convergence}
Let $S$ be a compact complex manifold with fundamental
group $\Gamma$, let $X$ be a compact $\Gamma$-space of
convergence type, and let $(M,\mathcal{F})$ be the
corresponding foliated bundle  over $S$.
Then $\mathcal{F}$ is holomorphically plain.
\end{propo}

A   class of $\Gamma$-spaces for which the convergence property
is well known to hold consists of actions
of non relatively compact subgroups of a Gromov-hyperbolic group
on the boundary of the latter. (See, for example, \cite{bow} and
references cited there. It should be noted that the standard
definition of the  convergence property used in the literature on
hyperbolic groups is much more restrictive than the one
we are using here.)  Therefore (keeping in mind Corollary
\ref{transmeaslemma}, the following holds.

\begin{propo}\label{hyperbolic}
Let $G$ be a Gromov-hyperbolic group,
  $X$   the boundary of $G$, and
  $S$   a compact connected complex   manifold with
fundamental group $\Gamma$. Suppose that $\Gamma$ acts on
$X$ via a homomorphism $\rho:\Gamma\rightarrow G$ and let
$(M_\rho,\mathcal{F}_\rho)$ be the corresponding foliated
bundle  over
$S$. Then $\mathcal{F}_\rho$ is holomorphically plain.
\end{propo}

\begin{coro}\label{rank1}
Let $S$ be a compact Riemann surface, let 
$\rho:\Gamma\rightarrow G$ be a homomorphism of
the fundamental group of $S$ into a connected simple Lie group
of rank one,
 and let
$(M_\rho,\mathcal{F}_\rho)$  be the foliated bundle over $S$
with fibers $X$, where $X$ is the boundary at infinity of 
the Riemannian symmetric space associated to $G$.
Then $\mathcal{F}_\rho$ is
holomorphically plain.
\end{coro}

\subsection{Codimension-One}
The main idea used in the proof of Proposition \ref{prop1},
together
with elementary facts about the structure of
codimension-one foliations yield the following. 

\begin{theor}\label{cod1}
If $(M,\mathcal{F})$ has codimension $1$, then it is
holomorphically plain.  
\end{theor}

\subsection{An Example}
The   results described so far might lead one to expect
 that holomorphically foliated  spaces are
holomorphically plain
 under very general conditions
and that one should be able to prove it
 using only   qualitative
properties of leafwise holomorphic functions.
 The next theorem shows, however, that the situation cannot
be so simple.

\begin{theor}\label{counterexample}
There exists a compact real   analytic foliation
$(M,\mathcal{F})$, which is a foliated bundle over
a compact Riemann surface,
 and a real analytic leafwise holomorphic function  
$f:M\rightarrow \mathbb{C}$ that  is not leafwise constant.
\end{theor}

To construct an example of $(M,\mathcal{F})$ and $f$ as in
Theorem \ref{counterexample} we first introduce some notation.
Let $\mathbb{D}$ denote as before the unit open disk in
$\mathbb{C}$. Points in  projective space $\mathbb{R}P^4$
will be written
$[z_1,z_2,
t]$, where  $z_i\in \mathbb{C}, t\in \mathbb{R}$ and  $(z_1,z_2,t)$ is
non-zero.
Define a (real analytic) action of $SU(1,1)$  on
$\mathbb{R}P^4$ as follows. 
Elements of $SU(1,1)$ are matrices of the form
$\begin{pmatrix}\alpha&\beta\\
\bar{\beta}&\bar{\alpha}
\end{pmatrix}$ for which $|\alpha|^2-|\beta|^2=1$.
The action of $SU(1,1)$ on $\mathbb{R}P^4$  defined by
 \begin{equation*}
 \begin{pmatrix}\alpha&\beta\\
\bar{\beta}&\bar{\alpha}
\end{pmatrix}\cdot [z_1,z_2,t]:=[\alpha z_1 +\beta \bar{z}_2, \alpha z_2
+\beta\bar{z}_1,t] 
\end{equation*}
leaves invariant the submanifold
$ \mathcal{C}:=\{[z_1,z_2,t]\in \mathbb{R}P^4:
|z_1|^2-|z_2|^2=t^2\}.$

We define on $\mathbb{D}\times\mathcal{C}$ the function
\begin{equation*}f(z,[\alpha, \beta, t]):=\frac{\bar{\alpha}
z-\beta}{-\bar{\beta}z +{\alpha}}.
\end{equation*}
An elementary calculation shows that $f(gz, g[\alpha,\beta, t])=f(z,
[\alpha,\beta, t])$ for every $g\in SU(1,1)$. Therefore,
if $\Gamma$ is a uniform lattice in $SU(1,1)$, then
$f$ yields a function on the foliated bundle $M=(\mathbb{D}\times
\mathcal{C})/\Gamma$ that is real analytic, leafwise holomorphic
and the restriction of $f$ to any leaf for which $t\neq 0$ is not
constant.

\subsection{A Universal Non-Plain Foliated Space}
We describe now a kind of ``universal space'' from which
such examples can be constructed. 
This will be done in the setting of  foliated bundles whose base are
compact Riemann surfaces, but it should be apparent that the same ideas
apply more broadly.
Let
$X_0:=\text{Hol}(\mathbb{D},\overline{\mathbb{D}})$ be the
space of holomorphic functions
 defined on $\mathbb{D}$ 
  such that $\sup\{|f(z)| : z\in
\mathbb{D}\}\leq 1.$ Then $X_0$, with the topology of
uniform convergence on compact sets, is a compact 
metrizable space upon which $PSU(1,1)$ acts via the
  continuous action: $(g,f)\mapsto    f\circ g^{-1}$,
where $g$, on the right-hand side, is regarded as an automorphism
(a M\"obius transformation) of
the Poincar\'e disc.

Let now $\rho:\Gamma\rightarrow PSU(1,1)$ be  a  
homomorphism from the fundamental group of a compact Riemann
surface $S=  \mathbb{D}/\Gamma$ into the M\"obius group,
and construct the   foliated bundle 
$(\tilde{S}\times X_0)/\Gamma$
over $S$.   We will denote
the resulting foliated space by $(M_0,\mathcal{F}_0)$.
By a {\em morphism} $f:(M,\mathcal{F})\rightarrow
(M_0,\mathcal{F}_0)$ we will mean a
(continuous)
$f:M\rightarrow M_0$ that maps leaves to leaves
holomorphically such that
$\pi_{M_0}\circ f=\pi_{M}$, where $\pi_M$ (resp.,
$\pi_{M_0}$) is the natural projection from $M$ to $S$
(resp.,
from $M_0$ to $S$).

A leafwise nonconstant, leafwise holomorphic
function can now be produced
on $(M_0,\mathcal{F}_0)$ by the following
essentially tautological procedure. First define
$\bar{\phi}:\mathbb{D}\times X_0\rightarrow \mathbb{C}$
by $\bar{\phi}(z,f):=f(z)$. Note that 
 $\bar{\phi}( \gamma(z),f\circ \gamma^{-1})=\bar{\phi}(z,f)$
for each $\gamma\in \Gamma$. 
There is as a result a well-defined function 
$\phi:M_0\rightarrow \mathbb{C}$ such that $\phi\circ \pi=\bar{\phi}$,
where $\pi$ is the natural projection from $\mathbb{D}\times
X\rightarrow M_0$.
The function $\phi:M_0\rightarrow \mathbb{C}$ is  a
continuous, leafwise holomorphic function. 

The following remark is an immediate 
consequence of these definitions.
In the proposition, equivariance of a map
$\hat{\psi}:V\rightarrow \mathbb{C}$
 means that
$\hat{\psi}\circ
\gamma=
\gamma\circ \hat{\psi}$ for each $\gamma\in \Gamma$.
\begin{propo}
 Let  $(M,\mathcal{F})$ be a foliated
bundle over $\Gamma\backslash \mathbb{D}$  with fiber $V$.
Then 
  there is a one-to-one correspondence between
(continuous) leafwise holomorphic functions  
$\psi:M\rightarrow
\mathbb{C}$  and  $\Gamma$-equivariant (continuous) 
$\hat{\psi}:V\rightarrow X_0$. Furthermore, if
${\Psi}:M\rightarrow M_0$
 is the morphism of holomorphically foliated spaces
induced from $\hat{\psi}$, then $\psi=\phi\circ
{\Psi}$, and $\Psi$ is the unique morphism from
$M$ to $M_0$ that satisfies this last equality.
\end{propo}
 
The   proposition indicates  how to go about looking for
examples of foliated manifolds that are not holomorphically plain: one
tries to find a
$\Gamma$-invariant manifold $V$ embedded in $X$. Specifically,
one can try to obtain a manifold $V\subset X$   as the closure
of a
$PSU(1,1)$-orbit. 
In fact, the example given just after Theorem \ref{counterexample}
is closely related  to    
what one gets by taking $V$ to be the closure of
the $PSU(1,1)$-orbit of the function $\varphi(z)=z$ in
$\text{Hol}(\mathbb{D},\overline{\mathbb{D}})$. 
This closure is the compactification of a $3$-manifold
(an open solid torus) by adding a circle at infinity (so as to
form a $3$-sphere), while the submanifold
$\mathcal{C}\subset\mathbb{R}P^4$  of that example
corresponds to  an analytic ``doubling'' of this
$3$-manifold.

The   remark just made suggests that a precise
characterization of holomorphically plain foliated bundles
over a compact Riemann surface will require an
investigation of 
  the dynamics of the   action 
of
$PSU(1,1)$ 
on $\text{Hol}(\mathbb{D},\overline{\mathbb{D}})$.
Such a characterization should tell, in particular,
how common holomorphically  plain foliations are,
at least in this special setting.

\section{Proofs}
\subsection{Proposition \ref{easy} and Corollary \ref{coroeasy} }
The hypothesis that $T\mathcal{F}$ is holomorphically
trivial means the following: there exist   vector 
   fields,  
$X_1,\dots,X_l$, on $M$, where $l$ is the leaf dimension,
such that the $X_i$ are everywhere tangent to
$\mathcal{F}$, linearly independent, and define  
holomorphic vector fields on leaves. Furthermore, 
the $X_i$, together with their tangential derivatives of
first order, are continuous on $M$.

Since $M$ is compact, the $X_i$ are complete vector fields
and each flow line is the image of a holomorphic map from
$\mathbb{C}$ into a leaf. Therefore, the restriction to
orbits of leafwise holomorphic functions define
bounded holomorphic functions on $\mathbb{C}$. As a result,
such functions
are constant on orbits, hence leafwise constant.

To show the corollary,
 write $[X_i,X_j]=\sum_k
f_{ij}^kX_k$. The coefficients
$f_{ij}^k$ are
continuous, leafwise holomorphic, hence leafwise constant.
   Due to the existence of a dense leaf, these
coefficients are constant on $M$. Therefore, the $X_i$
span a  finite dimensional Lie algebra.
The corollary is now a result of standard facts in Lie
theory.

\subsection{Proposition \ref{transmeas}}
A leafwise holomorphic function $f$ is
$\Delta_{\bar{\partial}}$-harmonic since $\bar{\partial}f=0$,
so the proposition is an immediate result of
\cite[Theorem 1(b)]{garnett}. We give here an independent
elementary proof when $(M,\mathcal{F})$ is leafwise K\"ahler.
It is well known (see, for example, \cite[p. 115]{griffiths})
that under this extra hypothesis
$\Delta_{\bar{\partial}}\bar{f}=
\Delta_{{\partial}}\bar{f}=0$, so we also have
$\Delta_{\bar{\partial}}\bar{f}=0$.
Another elementary 
calculation gives the identity  
$$\Delta_{\bar{\partial}}(|f|^2)=-\|\partial f\|^2 
 +f\Delta_{\bar{\partial}}\bar{f}=-\|\partial
f\|^2$$ holds.  By integrating
the resulting equality against 
  a harmonic measure $m$ one deduces that $\partial f=0$
  $m$-almost everywhere, and as $\partial f$ is continuous,
it is zero on the support of $m$. Therefore,
$f$ must be constant on that support.
 \subsection{Proposition \ref{prop1}}
The next lemma and corollary will be needed a number of
times in the rest of the paper.
 \begin{lem}\label{keylemma}
Let $(M,\mathcal{F})$ be a holomorphically foliated
compact connected space and let  $f:M\rightarrow
\mathbb{C}$ be a leafwise holomorphic function. Let
$C\subset M$ be the set where $f$ is leafwise constant.
Suppose that
$f(C)\subset
\mathbb{C}$ is at most countable. Then $f$ is constant
on $M$.
\end{lem}
\begin{pf}
Clearly $C$ is a compact $\mathcal{F}$-saturated set.
By the open mapping theorem for holomorphic functions,
the restriction of $f$ to each leaf in the complement
of $C$ is open, hence $f$ itself is an open mapping
on that complement. Therefore
$U:=f(M\backslash C)$ is a bounded open subset of
$\mathbb{C}$. Since $M$ is compact, $f(M)=U\cup f(C)$
is compact. In particular, the boundary of $U$ is contained
in $f(C)$, which is by assumption a countable set. But
a bounded non-empty open set in $\mathbb{C}$ cannot
have a countable boundary. Therefore $U$ is empty
and (as $M$ is connected) $f(M)$ reduces to a point.
\end{pf}

\begin{coro}\label{keycorollary}
Suppose that $(M,\mathcal{F})$ has (at most) countably 
many minimal sets. Then any leafwise holomorphic function 
is constant on $M$.
\end{coro}
\begin{pf}
Let $f$ be a leafwise holomorphic function and
$C$ as in Lemma \ref{keylemma}.
Denote by $\mathcal{M}$ the union of minimal sets.
We claim that $f(C)=f(\mathcal{M})$. Indeed if $x\in C$
and $f(x)=c$, then $f$ takes the constant value $c$
on the closure of the leaf containing $x$. But this
closure contains a minimal set, so $f(C)\subset
f(\mathcal{M})$. Conversely,
by an application of the maximal principle for holomorphic
functions, each closed saturated set
 contains leaves where
$f$ is constant, so $f(C)=f(\mathcal{M})$ as claimed.
The main assertion is now a consequence of Lemma
\ref{keylemma}.
\end{pf}

By applying Corollary \ref{keycorollary} to
the closure of each leaf we obtain  Proposition
\ref{prop1}. 

\subsection{Proposition \ref{generic}}
 
By \cite{rapin}, $\text{Hom}(\Gamma, G)$
is irreducible,  and by
\cite[8.2]{burger} the homomorphisms with Zariski
dense image form a Zariski open subset  
$\mathcal{D}\subset\text{Hom}(\Gamma, G)$. On the other
hand, $\mathcal{D}$ is clearly nonempty. ($\Gamma$ has
homomorphisms onto the free group on $2$ generators
and it is easy to show that the free group has 
representations with Zariski dense image.)
Therefore, the theorem is a 
   consequence of
Corollary   \ref{zariskidense}.

\subsection{Corollary \ref{convergence}}
We identify $X$ with a fixed fiber of the foliated bundle,
so that the holonomy transformations of $\mathcal{F}$ correspond
to the $\Gamma$-action on $X$.

Suppose that $f$ is a leafwise holomorphic continuous function
on $M$.
Let $C\subset X$ be the compact $\Gamma$-invariant subset
that corresponds to leaves on which $f$ is constant. Then 
there is a countable set of points $x_1, x_2, \dots$
such that for
all $y\in C$, with the possible exception of a countable
subset of $C$, we can find a sequence $\gamma_m$ 
such that $\gamma_mx\rightarrow x_i$ for some $i$.
Clearly the $x_i$ belong to $C$. Therefore $f$ can
take at most a countable set of  values on leaves of
$\mathcal{F}$ in $C$. So $f$ must be leafwise constant
by Lemma \ref{keylemma}.

\subsection{Theorem \ref{cod1}}
The basic facts about codimension-one foliations
that we will use can be found in \cite{hectorhirsch}
or \cite{candelconlon}, for example. We recall
that a minimal set is said
to be {\em exceptional} if
it is neither a single closed leaf of (a codimension-one
foliation) $\mathcal{F}$ nor a connected component
of $M$.

By a theorem of Haefliger \cite{haefliger} (see also
\cite{hectorhirsch}), the
union of  compact leaves of $\mathcal{F}$ is a compact
set, which we denote by $N$. Let $L$ be a connected
component    of the complement of $N$. Then, as
$\mathcal{F}$ has codimension $1$,  the closure
$\bar{L}$
 is a compact manifold whose boundary is a finite union
of compact leaves (cf. \cite{hectorhirsch}).

It is known that $(M,\mathcal{F})$ has only a finite
set of exceptional minimal sets. In particular,
$(\bar{L},\mathcal{F}|_{\bar{L}})$  has finitely
many minimal sets (exceptional or not).
Applying Corollary \ref{keycorollary} we deduce that any
leafwise holomorphic function is constant on $\bar{L}$.
In particular, any leafwise holomorphic function is
leafwise constant on $M\backslash N$. The same is obviously
true on $N$ (which is the union of closed leaves).
Therefore $(M,\mathcal{F})$ is holomorphically plain.


\begin{thebibliography}{Klg9}
\bibitem{burger} N. A'Campo and M. Burger. {\em R\'eseaux
arithm\'etiques et commensurateur d'apr\`es G. A.
Margulis},
Invent. math. {\bf 116}, 1-25 (1994)
\bibitem{candel} A. Candel. {\em The Harmonic measures of Lucy
Garnett}, preprint, 2000.
\bibitem{candelconlon} A. Candel and L. Conlon.
{Foliations I}, Graduate Studies in Mathematics, Volume 23,
AMS, 2000.
\bibitem{ghys2} D. Cerveau, E. Ghys, N. Sibony, J.C. Yoccoz.
Dynamique et G\'eom\'etrie Complexes, Panoramas et Synth\`eses,
Soci\'et\'e Math\'ematique de France, 1999.
\bibitem{connes} A. Connes. {\em A survey of foliations
and operator algebras}, Proc. Symp. Pure Math., Amer.
Math. Soc.,  521-628 (1982).
\bibitem{garnett} L. Garnett. {\em Foliations, the
ergodic theorem and Brownian motion}, J. Funct. Anal. {\bf
51} (1983), 285-311.
\bibitem{bow} E. Ghys and P. de la Harpe. Sur les groups
hyperboliques d'apr\`es Mikhael Gromov, Progress in Mathematics
83, Birkh\"auser, Basel, 1990. 
\bibitem{griffiths} P. Griffiths and J. Harris.
{Principles of Algebraic Geometry}, John Wiley \& Sons,
1994.
\bibitem{haefliger} A. Haefliger. {\em Vari\'et\'es
feuillet\'ees}, Ann. Scuola Norm. Sup. Pisa {\bf 16}
(1962),367-397.
\bibitem{hectorhirsch} G. Hector and U. Hirsch.
{Introduction to the geometric theory of foliations},
Aspects of Mathematics, 1983, Vieweg.
\bibitem{kodaira} K. Kodaira. {Complex manifolds and
deformations of complex structures}, Springer, 1986.
\bibitem{margulis} G. A. Margulis. {Discrete Subgroups of
Semisimple Lie Groups}, Springer, 1989.
\bibitem{molino} Pierre Molino. {Riemannian Foliations},
Birkhauser, 1987.
\bibitem{rapin} A. S. Rapinchuk, V. V. Benyash-Krivetz, V.
I. Chernousov. {\em Representation varieties of the
fundamental groups of compact orientable surfaces}, Israel
J. Math. {\bf 93} (1996) 29-71.
\end{thebibliography}
\end{document}